\documentclass[12pt]{article}

\usepackage{amsmath,amsthm,amsfonts,amssymb}

\newtheorem{theorem}{Theorem}[section]

\newtheorem{lemma}[theorem]{Lemma}
\newtheorem{definition}[theorem]{Definition}
\newtheorem{remark}[theorem]{Remark}

\def\cB{\mathcal{B}}
\def\cE{\mathcal{E}}
\def\cF{\mathcal{F}}
\def\cP{\mathcal{P}}

\def\bD{\mathbb{D}}
\def\bR{\mathbb{R}}
\def\bS{\mathbb{S}}
\def\bZ{\mathbb{Z}}

\def\barR{\overline{\bR}}

\def\barD{\overline{\bD}}

\begin{document}

\title{Regular variation of infinite series of processes with random coefficients}

\author{Raluca Balan\footnote{Department of Mathematics and Statistics, University of Ottawa,
585 King Edward Avenue, Ottawa, ON, K1N 6N5, Canada. E-mail
address: rbalan@uottawa.ca} \footnote{Research supported by a
grant from the Natural Sciences and Engineering Research Council
of Canada.}}

\date{January 30, 2014}
\maketitle

\begin{abstract}
\noindent In this article, we consider a series $X(t)=\sum_{j \geq 1}\Psi_j(t) Z_j(t),t \in [0,1]$ of random processes with sample paths in the space $\bD$ of c\`adl\`ag functions (i.e. right-continuous functions with left limits) on $[0,1]$. We assume that $(Z_j)_{j \geq 1}$ are i.i.d. processes with sample paths in $\bD$ and $(\Psi_j)_{j \geq 1}$ are processes with continuous sample paths. Using the notion of regular variation for $\bD$-valued random elements (introduced in \cite{hult-lindskog05}), we show that $X$ is regularly varying if $Z_1$ is regularly varying, $(\Psi_j)_{j \geq 1}$ satisfy some moment conditions, and a certain ``predictability assumption'' holds for the sequence $\{(Z_j,\Psi_j)\}_{j \geq 1}$.
Our result can be viewed as an extension of Theorem 3.1 of \cite{hult-samorod08} from random vectors in $\bR^d$ to random elements in $\bD$. As a preliminary result, we prove a version of Breiman's lemma for $\bD$-valued random elements, which can be of independent interest.
\end{abstract}

\noindent {\em MSC 2010:} Primary 60G52; secondary 60G17, 62M10


\vspace{3mm}

\noindent {\em Keywords:} regular variation, linear processes, stochastic recurrence equation, Breiman's lemma

\section{Introduction}

Regular variation is an important property which lies at the core of some fundamental results in probability theory, which describe the asymptotic behavior of the maximum and the sum of $n$ i.i.d. random variables. In the past 30 years, especially after the publication of the landmark articles \cite{davis-resnick85,resnick86} and monograph \cite{resnick87}, this has become a very active area of research, with a huge potential for applications, arising usually in the context of time series models.

A random variable $Z$ is {\em regularly varying} if $P(|Z|>x)=x^{-\alpha}L(x)$ for any $x>0$, and $P(Z>x)/P(|Z|>x) \to p$ as $x \to \infty$, for some $\alpha>0,p \in [0,1]$ and a slowly varying function $L$. This is equivalent to the fact that $|Z|$ is in the maximal domain of attraction of the Fr\'echet distribution $\Phi_{\alpha}(x)=\exp(-x^{-\alpha})$, and if $\alpha<2$, to
the fact that $Z$ is in the domain of attraction of a stable distribution with index $\alpha$. Another useful characterization states that $Z$ is regularly varying if and only if there exists a sequence $(a_n)_n \uparrow \infty$ such that $nP(|Z|>a_n x) \to x^{-\alpha}$ for any $x>0$. By considering the state space $\barR_0=[-\infty,\infty] \verb2\2 \{0\}$ instead of $\bR$ (such that sets of the from $[-\infty,-x) \cup (x,\infty],x>0$ become relatively compact), the previous convergence can be expressed as the vague convergence of Radon measures:
$$nP(a_n^{-1}Z \in \cdot) \stackrel{v}{\to}\nu(\cdot) \quad \mbox{in} \ \barR_0,$$ where
$\nu(dx)=[p\alpha x^{-\alpha-1}1_{(0,\infty)}(x)+(1-p)\alpha (-x)^{-\alpha-1}1_{(-\infty,0)}(x)]dx$
is a measure on $\barR_0$ with $\nu(\barR_0 \verb2\2 \bR)=0$ (see e.g. Section 3.6 of \cite{resnick07}).

Recall that a measure $\mu$ on a locally compact space with a countable basis (LCCB) is {\em Radon} if $\mu(B)<\infty$ for any relatively compact Borel set $B$.
A sequence $(\mu_n)_n$ of Radon measures {\em converges vaguely} to a Radon measure $\mu$ (written as $\mu_n \stackrel{v}{\to} \mu$) if $\mu_n(B) \to \mu(B)$ for any relatively compact Borel set $B$ with $\mu(\partial B)=0$ (see Appendix 15.7 of \cite{kallenberg83}).

In higher dimensions, a random vector ${\bf Z}$ with values in $\bR^d$ is called {\em regularly varying} if there exist a sequence $(a_n)_n \uparrow \infty$ and a non-null Radon measure $\nu$ on $\barR_0^d=[-\infty,\infty]^d \verb2\2 \{{\bf 0}\}$ such that $\nu(\barR_0^d \verb2\2 \bR^d)=0$ and
\begin{equation}
\label{RV}
nP(a_n^{-1}{\bf Z} \in \cdot) \stackrel{v}{\to} \nu(\cdot) \quad \mbox{in} \quad \barR_0^d.
\end{equation}
In this case, we write $Z \in {\rm RV}(\{a_n\},\nu, \barR_0^d)$.
It can be proved that the measure $\nu$ satisfies the following {\em scaling property}: there exists $\alpha>0$ (called the index of ${\bf Z}$) such that
$\nu(sB)=s^{-\alpha} \nu(B)$ for any $s>0$ and for any Borel set $B \subset \barR_0^d$.
In particular, for any norm $\|\cdot\|$ on $\bR^d$ and for any $r>0$,
$$n P(\|{\bf Z}\|>a_n r) \to cr^{-\alpha},$$
where $c=\nu(\{{\bf x} \in \barR_0^d;\|{\bf x}\| >1\})$.
Let $\bR_0^d=\bR^d \verb2\2\{{\bf 0}\}$, $\bS=\{{\bf x} \in \bR^d; \|{\bf x}\|=1\}$ be the unit sphere in $\bR^d$ and $T:\bR_0^d  \to (0,\infty) \times \bS$ be the polar coordinate transformation: $T({\bf x})=(\|{\bf x}\|,{\bf x}/\|{\bf x}\|)$.
The fact that ${\bf Z} \in {\rm RV}(\{a_n\},\nu,\barR_0^d)$ can also be expressed as: (see e.g. Section 6.1 of \cite{resnick07})
\begin{equation}
\label{RVpolar}
nP(T(a_n^{-1}{\bf Z}) \in \cdot) \stackrel{v}{\to} c\nu_{\alpha} \times \sigma \quad \mbox{in} \quad (0,\infty] \times \bS
\end{equation}
where $\nu_{\alpha}(r,\infty)=r^{-\alpha}$ and $\sigma$ is a probability measure on $\bS$ given by $\sigma(S)=c^{-1}\nu(\{{\bf x} \in \barR_0^d; \|{\bf x}\|>1,{\bf x}/\|{\bf x}\| \in S\})$.
Hence, $\nu \circ T^{-1}=c\nu_{\alpha} \times \sigma$ on $(0,\infty) \times \bS$.

In the one dimensional case, many time series models can be expressed as linear processes of the form:
\begin{equation}
\label{series}
X_i=\sum_{j \geq 0}c_j Z_{i-j}, \quad i \in \bZ
\end{equation}
where $(c_j)_{j \geq 1}$ are real numbers and $(Z_{j})_{j \in \bZ}$ are i.i.d. random variables. One simple example is the auto-regressive model of order 1, $X_i=a X_{i-1}+Z_i$ with $|a|<1$, leading to $X_i=\sum_{j \geq 0}a^{j}Z_{i-j}$.
Assume that $Z_0$ is
regularly varying with index $\alpha$ and slowly varying function $L$ as above. A basic question is: under what conditions the series \eqref{series} converges and if so, is $X_0$ still regularly varying? If $\alpha<2$, an argument which can be traced back to \cite{astrauskas83} (see also Proposition 2.1 of \cite{BJL}) shows that the series \eqref{series} converges a.s. if and only if
$$\sum_{j \geq 0}|c_j|^{\alpha}L(1/|c_j|)<\infty.$$
In \cite{mikosch-samorod00}, it was shown that sufficient conditions for the converges of the series \eqref{series} are:
$\sum_{j \geq 1}|c_j|^{\alpha-\gamma}<\infty$  for some $\gamma \in (0,\alpha)$, if $\alpha \leq 2$, and $\sum_{j \geq 0}|c_j|^2<\infty$ if $\alpha>2$, and
under these conditions, $X_0$ is regularly varying. This result continues to hold for $d$-dimensional vectors $(Z_j)_j$, and deterministic $p \times d$ matrices $(A_j)_j$ replacing the coefficients $(c_j)_j$  (see Corollary 3.1 of \cite{hult-samorod08}).

More interesting models lead to series of the form:
\begin{equation}
\label{series2}
X_{i}=\sum_{j \geq 0}C_{i,j}Z_{i-j}, \quad i \in \bZ
\end{equation}
with random coefficients $C_{i,j}$. One such example is the stochastic recurrence equation (SRE)
$X_i=Y_i X_{i-1} +Z_{i}$, where $\{(Y_i,Z_i)\}_{i \in \bZ}$ is a sequence of i.i.d. random vectors in $\bR^2$.  Under certain conditions, it can be shown that (SRE) has a unique stationary solution which can be represented as a series of the form \eqref{series2} with $C_{i,0}=1$ and $C_{i,j}=\prod_{k=1}^{j}Y_{i-k+1}$ for $j \geq 1$. Another example is the bilinear model $X_i=cX_{i-1}Z_{i-1}+Z_i$, which admits a stationary solution of the form $X_i=Z_i+\sum_{j \geq 1}C_{i,j}Z_{i-j}^2$ with $C_{i,j}=\prod_{k=1}^{j-1}Z_{i-k}$.
Assume that $Z_0$ is regularly varying. The first result addressing the question mentioned above is due to \cite{resnick-willekens91}, where it was shown that a series $X_0$ given by \eqref{series2} is regularly varying, if $(C_{0,j})_{j \geq 0}$ and $(Z_{-j})_{j \geq 0}$ are independent and $(C_{0,j})_{j \geq 0}$ satisfy some moment conditions. This result was extended to (SRE) and the bilinear model in \cite{grey94}, respectively \cite{davis-resnick96}. This was improved in \cite{hult-samorod08}, under weaker moment conditions on $(C_{0,j})_{j \geq 0}$ and a certain ``predictability assumption'' imposed on the sequence $\{(C_{0,j},Z_{-j})_{j \geq 0}$. The result of \cite{hult-samorod08} is in fact valid for random vectors, one of the major applications of (SRE) in higher dimensions being the GARCH model (see \cite{basrak-et-al02}).

A breakthrough idea, which gave a new perspective to the concept of regular variation and lead to a different line of investigations, was introduced in \cite{dehaan-lin01}. Motivated by extreme value theory, this idea was to examine the {\em global} asymptotic behavior (as $t$ runs in a fixed interval $[0,1]$) of the normalized maximum process $\{a_n^{-1}\max_{1 \leq i \leq n}Z_i(t)\}_{t \in [0,1]}$ associated with $n$ i.i.d. processes $Z_1, \ldots,Z_n$ whose sample paths are c\`adl\`ag functions (i.e. continuous functions with left limits) on $[0,1]$. Each process $Z_i$ is interpreted as a collection of measurements observed continuously over a fixed linear spatial region, $Z_i(t)$ being the observation at location $t$ and time $i$. In the example of \cite{dehaan-lin01}, $Z_i(t)$ is the high tide water level at location $t$ and time $i$, along the northern coast of the Netherlands.
It turns out that if $Z_0$ is regularly varying (in a sense which is made precise in Section \ref{regvarD-section} below), then the finite-dimensional distributions of the normalized maximum process converge to those of a c\`adl\`ag process $Y=\{Y(t)\}_{t \in [0,1]}$, and if $Y$ has continuous sample paths, then the convergence is in distribution in the space
$\bD=\bD[0,1]$ of c\`adl\`ag functions on $[0,1]$, endowed with the Skorohod $J_1$-topology
(see Theorem 2.4 of \cite{dehaan-lin01}).

The notion of regular variation for c\`adl\`ag processes was thoroughly studied in \cite{hult-lindskog05} where it was proved that it is equivalent to the regular variation of the finite-dimensional distributions of the process, combined with some relative compactness conditions
(see Theorem 10 of \cite{hult-lindskog05}). In particular, a L\'evy process $\{Z(t)\}_{t \in [0,1]}$ is regularly varying if and only if $Z(1)$ is regularly varying (see Lemma 2.1 of \cite{hult-lindskog07}). One example is the $\alpha$-stable L\'evy process.

In this context, it becomes interesting to examine the regular variation of $\bD$-valued time series of the form:
\begin{equation}
\label{series3}
X_i(t)=\sum_{j \geq 0}\psi_j(t)Z_{i-j}(t), \quad t \in [0,1],i \in \bZ
\end{equation}
where $Z_i=\{Z_i(t)\}_{t \in [0,1]}, i \in \bZ$ are i.i.d. regularly varying processes and $\psi_j=\{\psi_j(t)\}_{t \in [0,1]},j \geq 0$ are deterministic functions in $\bD$. This analysis was carried out in \cite{davis-mikosch08}, where it was proved that, under some conditions on the coefficients $(\psi_j)_{j \geq 0}$, the series $X_0$ given by \eqref{series3} converges a.s. and is regularly varying in $\bD$.
Moreover, the authors of \cite{davis-mikosch08} derived the limit distribution of the normalized space-time maximum $a_n^{-1}\max_{i \leq n}\sup_{t \in [0,1]}|X_i(t)|$, using non-trivial point process techniques.
We should note that the results of \cite{davis-mikosch08} are in fact valid for c\`adl\`ag processes indexed by $[0,1]^d$ with $d \geq 1$, being motivated by applications to spatial processes.

In the present article, we consider the next natural step in this line of investigations which consists in examining series of the form:
\begin{equation}
\label{series4}
X_{i}(t)=\sum_{j \geq 0}\Psi_{i,j}(t)Z_{i-j}(t), \quad t \in [0,1],i \in \bZ
\end{equation}
where  $Z_i=\{Z_i(t)\}_{t \in [0,1]}, i \in \bZ$ are i.i.d. regularly varying processes and $\Psi_{i,j}=\{\Psi_{i,j}(t)\}_{t \in [0,1]}$ are random processes. For example, one can consider that at each spatial location $t$, the temporal dependence between the observations is described by an (SRE) model $X_i(t)=Y_i(t)X_{i-1}(t)+Z_i(t)$, leading to model \eqref{series4} with $\Psi_{i,0}(t)=1$ and $\Psi_{i,j}(t)=\prod_{k=1}^{j}Y_{i-k+1}(t)$ for $j \geq 1$.
Our main result shows that, if $(\Psi_{0,j})_{j \geq 0}$ satisfy some moment conditions, and the same ``predictability assumption'' as in \cite{hult-samorod08} holds for the sequence $\{(\Psi_{0,j},Z_{-j})\}_{j \geq 0}$, then the series $X_0$ given by \eqref{series4} converges a.s. and is regularly varying in $\bD$.
We postpone the asymptotic analysis of the normalized maximum of $X_1,\ldots,X_n$ for a future study.

The article is organized as follows. In Section \ref{regvarD-section}, we introduce the concept of regular variation on $\bD$ and discuss some of its properties. In Section \ref{main-section}, we state and prove our result about the regular variation in $\bD$ of a series of the form \eqref{series4}. For this, we use two preliminary results, one of them being a version of
Breiman's lemma for $\bD$-valued random elements. The appendix contains a variant of Pratt's lemma (regarding the interchanging of $\limsup$ with an integral), which is needed for checking the relative compactness conditions mentioned above.

\section{Regular variation on $\bD$}
\label{regvarD-section}

In this section, we recall the definition and main properties of the regular variation for random processes with sample paths in $\bD$. We follow references  \cite{lindskog04,hult-lindskog05,davis-mikosch08}.

We let $\bD=\bD[0,1]$ be the space of right continuous functions $x:[0,1] \to \bR$ with left limits. Recall that $\bD$ is a complete separable metric space (CSMS), equipped with a distance called $d_0$, which is equivalent to Skorohod $J_1$-metric (see pages 109-115 of \cite{billingsley68}). We denote by $\cB(\bD)$ the class of Borel sets in $\bD$, equipped with the $J_1$-topology. Note that $\|x\|_{\infty}=\sup_{t \in [0,1]}|x(t)|<\infty$ for any $x \in \bD$, and the topology of uniform convergence on $\bD$ is stronger than the $J_1$-topology.

We let $\bS_{\bD}=\{x \in \bD; \|x\|_{\infty}=1\}$ be the ``unit sphere'' in $\bD$, endowed with metric $d_0$, and $\cB(\bS_{\bD})$ be the class of Borel sets in $\bS_{\bD}$.
We denote $\bD_0=\bD \verb2\2 \{0\}$, where $0$ is the null function in $\bD$, and let $\cB(\bD_0)$ be the class of Borel sets in $\bD_0$. Similarly to the polar coordinate transformation in $\bR^d$, we consider the homeomorphism $T: \bD_0 \to (0,\infty) \times \bS_{\bD}$ given by $T(x)=(\|x\|_{\infty},x/\|x\|_{\infty})$.
We define the space $$\barD_0:=(0,\infty] \times \bS_{\bD}.$$ This space is endowed with the product metric, where $(0,\infty]$ has the metric $\rho(x,y)=(1/x)-(1/y)$, with the convention $1/\infty=0$.
We let $\cB(\barD_0)$ be the class of all Borel sets in $\barD_0$. Note that $\barD_0 \verb2\2 T(\bD_0)=\{\infty\} \times \bS_{\bD}$.

\begin{remark}
{\rm The authors of \cite{lindskog04,hult-lindskog05,davis-mikosch08} identify the space $\bD$ with $[0,\infty) \times \bS_{\bD}$ and write $\barD_0 \verb2\2 \bD=\{\infty\} \times \bS_{\bD}$. For the sake of the analogy with $\bR^d$, we prefer to distinguish between $\bD_0$ and $(0,\infty) \times \bS_{\bD}$.}
\end{remark}

Similarly to $\bR^d$, the concept of regular variation on $\bD$ can be defined using convergence of measures. A small technical issue is the fact that $\barD_0$ is not a LCCB space, and hence the notion of vague convergence is not appropriate on this space. Fortunately, $\barD_0$ is a CSMS and vague convergence can be replaced by the $\hat{w}$-convergence.
Recall that a measure $\mu$ on a CSMS $E$ (with metric $d$) is {\em boundedly finite} if $\mu(B)<\infty$ for any bounded Borel set $B$ in $E$. (A set $B$ is bounded if it is contained in an open sphere $S_{r}(x)=\{y \in E; d(x,y)<r\}$.)
A sequence $(\mu_n)_n$ of boundedly finite measures converges to a boundedly finite measure $\mu$ in the $\hat{w}$-topology (written as $\mu_n \stackrel{\hat{w}}{\to} \mu$) if $\mu_n(B) \to \mu(B)$ for any bounded Borel set $B$ with $\mu(\partial B)=0$ (see Appendix A2.6 of \cite{D-VJ}).

The following definition introduces the analogue of \eqref{RVpolar} for $\bD$. Let $\nu_{\alpha}$ be the measure on $(0,\infty]$ given by $\nu_{\alpha}(dx)=\alpha x^{-\alpha-1}1_{(0,\infty)}(x)dx$, $\nu_{\alpha}(\{\infty\})=0$.

\begin{definition}
\label{def-RV-D}
{\rm A process $Z=\{Z(t)\}_{t \in [0,1]}$ with sample paths in $\bD$ is called {\em regularly varying} if 
there exist $\alpha>0,c>0$, a sequence $(a_n)_{n \geq 1}$ with $a_n>0,a_n \uparrow \infty$, and a probability measure $\sigma$ on $\bS_{\bD}$ such that
$$nP(T(a_n^{-1} Z) \in \cdot) \stackrel{\hat{w}}{\to} c\nu_{\alpha} \times \sigma \quad \mbox{in} \quad \barD_0.$$
$\alpha$ is called the {\em index} of $Z$.
}
\end{definition}

\begin{remark}
{\rm The previous definition coincides with Definition 2 of \cite{hult-lindskog05}, if we identify $T(x)$ with $x$, for any $x \in \bD_0$. Note that $\mu=c\nu_{\alpha} \times \sigma$ is a non-null boundedly finite measure on $\barD_0$ which satisfies $\mu(\{\infty\} \times \bS_{\bD})=0$.
}
\end{remark}

\begin{remark}
{\rm We now examine the analogue of \eqref{RV} for $\bD$. Suppose that $Z$ is regularly varying as in Definition \ref{def-RV-D}. Let $\cP$ be the class of sets of the form
$$V_{a,b;S}=\{x \in \bD_0;a<\|x\|_{\infty} \leq b, x/\|x\|_{\infty} \in S \}=T^{-1}((a,b] \times S),$$
for some $0<a<b<\infty$ and $S \in \cB(\bS_{\bD})$. Define $\nu(V_{a,b;S})=c\nu_{\alpha}((a,b])\sigma(S)$.
Since $\cP$ is a semiring which generates $\cB(\bD_0)$, by Theorem 11.3 of \cite{billingsley95}, $\nu$ can be extended to a measure on $\bD_0$.
Let $V_{r;S}=\{x \in \bD;\|x\|_{\infty}>r, x/\|x\|_{\infty} \in S\}$ be the $\bD$-analogue of a ``pizza-slice'' set from $\bR^d$, with $0<r<\infty$. Then
\begin{equation}
\label{nu-of-V}
\nu(V_{r;S})=cr^{-\alpha}\sigma(S),
\end{equation}
and hence $c=\nu(V_{1; \bS_{\bD}})$.
It can be shown that $\nu$ satisfies the scaling property $\nu(sB)=s^{-\alpha} \nu(B)$ for any $s>0$ and $B \in \cB(\bD_0)$,
$\nu(\partial V_{r;S})=cr^{-\alpha}\sigma(\partial S)$ and
\begin{equation}
\label{conv}
nP(a_n^{-1}Z \in V_{r;S}) \to \nu(V_{r;S})
\end{equation}
for any $r>0$ and $S \in \cB(\bS_{\bD})$ with $\sigma(\partial S)=0$ (see the proofs of Theorems 1.14 and 1.15 of \cite{lindskog04} for $\bR^d$; the same arguments work for $\bD$).
But \eqref{conv} cannot be expressed as a statement of $\hat{w}$-convergence, because there is no natural ``infinity'' that can be added to $\bD_0$.
Taking $S=\bS_{\bD}$ in \eqref{conv}, we obtain that for any $r>0$,
\begin{equation}
\label{tail-conv}
nP(\|Z\|_{\infty}>a_n r) \to c r^{-\alpha}.
\end{equation}
By abuse of notation, we write $Z \in {\rm RV}(\{a_n\},\nu, \barD_0)$, although $\nu$ is a measure on $\bD_0$, not on $\barD_0$. 
Note that
$$\nu \circ T^{-1}=c\nu_{\alpha} \times \sigma \quad \mbox{on} \quad (0,\infty) \times \bS_{\bD}.$$
}
\end{remark}

\vspace{3mm}

To introduce another characterization of regular variation on $\bD$, we need to recall the definition of the modulus of continuity: for any $x \in \bD$ and $\delta>0$,
$$w(x,\delta)=\sup_{|s-t| \leq \delta}|x(s)-x(t)|.$$
If $x$ is continuous, then $\lim_{\delta \to 0}w(x,\delta)=0$.
In the case of a discontinuous function $x \in \bD$, the following quantity plays the same role as $w(x,\delta)$:
$$w''(x,\delta)=\sup_{t_1 \leq t \leq t_2, \ t_2-t_1 \leq \delta} |x(t)-x(t_1)| \wedge |x(t_2)-x(t)|,$$
since $\lim_{\delta \to 0}w''(x,\delta)=0$ (see (14.8) and (14.46) of \cite{billingsley68}). We define
$$w(x,T)=\sup_{s,t \in T}|x(s)-x(t)| \quad \mbox{for any set} \ T \subset [0,1].$$

The following result will be needed in the sequel. This result shows that the regular variation in $\bD$ coincides with the regular variation of the marginal distributions, combined with some relative compactness conditions.

\begin{theorem}[Theorem 10 of \cite{hult-lindskog05}]
\label{HL-theorem}
Let $Z=\{Z(t)\}_{t \in [0,1]}$ be a process with sample paths in $\bD$. The following statements are equivalent: \\
(i) $Z \in {\rm RV}(\{a_n\},\nu,\barD_0)$; \\
(ii) There exists a sequence $(a_n)_{n \geq 1}$ with $a_n>0,a_n \uparrow \infty$, a set $T \subset [0,1]$ containing $0,1$ with $[0,1] \verb2\2 T$ countable, and a collection $\{\nu_{t_1, \ldots,t_k}; t_1, \ldots,t_k \in T,k \geq 1\}$, each $\nu_{t_1,\ldots,t_k}$ being a Radon measure on $\barR_0^{k}$ with $\nu_{t_1, \ldots,t_k}(\barR_0^{k} \verb2\2 \bR_0^k)=0$ and $\nu_t$ is non-null for some $t \in T$, such that:\\
(a)
$(Z(t_1), \ldots,Z(t_k)) \in {\rm RV}(\{a_n\},\nu_{t_1, \ldots,t_k},\barR_0^k)$ for all $t_1, \ldots,t_k \in T$; and \\
(b) the following three conditions are satisfied:
\begin{eqnarray*}
& (C1) & \lim_{\delta \to 0} \limsup_{n \to \infty} n P(w''(Z,\delta)>a_n \varepsilon)=0 \ \mbox{for any} \ \varepsilon>0; \\
& (C2) & \lim_{\delta \to 0} \limsup_{n \to \infty} n P(w(Z,[0,\delta))>a_n \varepsilon)=0 \ \mbox{for any} \ \varepsilon>0; \\
& (C3) & \lim_{\delta \to 0} \limsup_{n \to \infty} n P(w(Z,[1-\delta,1))>a_n \varepsilon)=0 \ \mbox{for any} \ \varepsilon>0.
\end{eqnarray*}
\end{theorem}

\begin{remark}
{\rm The sequences $\{a_n\}$ in (i) and (ii) can be taken to be the same.
The measure $\nu$ is uniquely determined by $\{\nu_{t_1, \ldots,t_k};t_1, \ldots,t_k \in T,k \geq 1\}$ and
\begin{equation}
\label{fin-dim-nu}
\nu_{t_1, \ldots,t_k}(B)=\nu (\pi_{t_1,\ldots,t_k}^{-1}(B \cap \bR^k)), \quad \mbox{for all} \ B \in \cB(\barR_0^k)
\end{equation}
where $\pi_{t_1, \ldots,t_k}(x)=(x(t_1), \ldots,x(t_k)),x \in \bD$.
The set $T$ in {\em (ii)} can be taken to be the set of all $t \in [0,1]$ such that $\nu(\{x \in \bD_0; \mbox{$x$ is not continuous at $t$}\})=0$.
}
\end{remark}

\section{The main result}
\label{main-section}

We are now ready to state our main result.

\begin{theorem}
\label{RV-infinite-sum}
Let $Z_j=\{Z_j(t)\}_{t \in [0,1]},j \geq 1$ be i.i.d. processes with sample paths in $\bD$ such that $Z_1 \in {\rm RV}(\{a_n\},\nu,\barD_0)$, and $\alpha>0$ be the index of $Z_1$.
Let $\Psi_j=\{\Psi_j(t)\}_{t \in [0,1]},j \geq 1$ be some processes with continuous sample paths, such that $P(\cup_{j \geq 1} \{\|\Psi_j\|_{\infty}>0 \})=1$, and
there exists an $m \geq 1$ and a set $T_1 \subset [0,1]$ containing $0$ and $1$, with $[0,1]\verb2\2 T_1$ countable, for which
\begin{equation}
\label{nonzero-psi-sum}
P(\bigcup_{j=1}^{m}\{\Psi_j(t) \not=0\})>0 \quad \mbox{for all} \ t \in T_1.
\end{equation}
Suppose that $(\Psi_1,\ldots,\Psi_j,Z_1,\ldots,Z_{j-1})$ is independent of $(Z_{k})_{k \geq j}$
for any $j \geq 2$, $\Psi_1$ is independent of $(Z_{j})_{j \geq 1}$, and there exists $\gamma \in (0,\alpha)$ such that:
\begin{eqnarray*}
&  & \sum_{j =1}^{m}E\|\Psi_j\|_{\infty}^{\alpha-\gamma}<\infty \ \mbox{and} \ \sum_{j \geq 1}E\|\Psi_j\|_{\infty}^{\alpha+\gamma}<\infty \quad  \mbox{if} \quad \alpha \in (0,1) \cup (1,2),\\
& & E\left(\sum_{j \geq 1}\|\Psi_j\|_{\infty}^{\alpha-\gamma}\right)^{(\alpha+\gamma)/(\alpha-\gamma)}<
\infty
\quad  \mbox{if} \quad \alpha \in \{1,2\}, \ \mbox{or} \\
&  & E \left(\sum_{j \geq 1}\|\Psi_j\|_{\infty}^2 \right)^{(\alpha+\gamma)/2}<\infty
\quad  \mbox{if} \quad \alpha>2.
\end{eqnarray*}
Then the series $X=\sum_{j \geq 1}\Psi_j Z_j$ converges in $\bD$ a.s. and $X \in {\rm RV}(\{a_n\},\nu^{X},\barD_0)$ where $$\nu^X(\cdot)=E\left[\sum_{j \geq 1} \nu \circ h_{\Psi_j}^{-1}(\cdot) \right].$$
For any $\psi \in \bD$, we define the product map $h_{\psi}:\bD \to \bD$ by $h_{\psi}(x)=\psi x, x \in \bD$, with $(\psi x)(t)=\psi(t)x(t)$ for any $t \in [0,1]$.
\end{theorem}

We begin with some preliminary results. The following result is known in the literature as {\em Breiman's lemma} (see \cite{breiman65}).

\begin{lemma}
\label{breiman-lemma-R}
Let $Z$ and $Y$ be independent nonnegative random variables such that $Z \in {\rm RV}(\{a_n\},\nu,\barR_0)$ and
$0<E(Y^{\alpha+\gamma})<\infty$ for some $\gamma>0$, where $\alpha>0$ is the index of $Z$ (and hence, $\nu(r,\infty)=cr^{-\alpha}$ for any $r>0$, for some $c>0$). Then $X=YZ \in {\rm RV}(\{a_n\},\nu^X,\barR_0)$ where $\nu^{X}(r,\infty)=cr^{-\alpha}E(Y^{\alpha})$ for any $r>0$.
\end{lemma}

Note that in Breiman's lemma,
$\nu^X(\cdot)=E[\nu\circ h_{Y}^{-1} (\cdot \cap [0,\infty))]$,
 where $h_y(x)=yx$ for any $x,y \in [0,\infty)$.

Proposition A.1 of \cite{basrak-et-al02} gives an extension of Lemma \ref{breiman-lemma-R} to a product $X=AZ$, where $Z$ is regularly varying in $\bR^d$ (with index $\alpha$), $A$ is a random matrix with $0<E\|A\|^{\alpha+\gamma}<\infty$ for some $\gamma>0$, and $Z,A$ are independent. Lemma 4.3 of \cite{hult-samorod08} extends this result to a finite sum $X=\sum_{j=1}^{m}A_j Z_j$, where $(Z_j)_j$ are i.i.d. regularly varying in $\bR^d$, $(A_j)_j$ are random matrices with $E\|A_j\|^{\alpha+\gamma}<\infty$ for some $\gamma>0$, and $Z_j$ is independent of $(A_1, \ldots,A_j,Z_1, \ldots,Z_{j-1})$ for all $j$. (For the later result, one also needs
the hypothesis $P(\cup_{j=1}^{m}\{\|A_j\|>0\})>0$, which is missing from \cite{hult-samorod08}.)

Our first result is a version of Breiman's lemma for processes with sample paths in $\bD$.

\begin{lemma}
\label{product-lemma}
Let $Z=\{Z(t)\}_{t \in [0,1]}$ and $\Psi=\{\Psi(t)\}_{t \in [0,1]}$ be independent processes with sample paths in $\bD$ such that
$Z \in {\rm RV}(\{a_n\},\nu,\barD_0)$, $\Psi$ has continuous sample paths, and
$E\|\Psi\|_{\infty}^{\alpha+\gamma}<\infty$ for some $\gamma>0$,
where $\alpha>0$ is the index of $Z$.
Suppose that there exists a set $T_1 \subset [0,1]$ containing $0$ and $1$ with $[0,1] \verb2\2 T_1$ countable, such that
\begin{equation}
\label{nonzero-psi}
P(\Psi(t) \not= 0)>0 \quad \mbox{for all} \ t \in T_1.
\end{equation}
Then $X=\Psi Z \in {\rm RV}(\{a_n\},\nu^X,\barD_0)$ where
\begin{equation}
\label{def-nux-prod}
\nu^X(\cdot)=E[\nu \circ h_{\Psi}^{-1}(\cdot)].
\end{equation}
\end{lemma}

\noindent {\bf Proof:} Let $T \subset [0,1]$ and $\{\nu_{t_1, \ldots,t_k};t_1, \ldots,t_k \in T,k \geq 1\}$  be the set and the marginal measures given by Theorem \ref{HL-theorem}.(ii) for $Z_1$. We show that $X$ satisfies the conditions of Theorem \ref{HL-theorem}.(ii) with $T_X=T \cap T_1$ instead of $T$ and the measures $\nu_{t_1, \ldots,t_k}^X$ (defined by \eqref{marginals-nux-prod} below) instead of $\nu_{t_1,\ldots,t_k}$.

First we show that $X$ satisfies condition (a). For this, let $t_1, \ldots,t_k \in T_X$ be arbitrary. Note that $(X(t_1), \ldots,X(t_k))^T=AY$ where $A$ is the diagonal matrix with entries $\Psi(t_1), \ldots, \Psi(t_k)$ and $Y=(Z(t_1),\ldots,Z(t_k))^T$. By Proposition A.1 of \cite{basrak-et-al02}, $(X(t_1), \ldots,X(t_k)) \in {\rm RV}(\{a_n\},\nu_{t_1, \ldots,t_k}^{X},\barR_0^k)$,
where
\begin{equation}
\label{marginals-nux-prod}
\nu_{t_1, \ldots,t_k}^X(\cdot)=E[\nu_{t_1, \ldots,t_k} \circ h_A^{-1}(\cdot \cap \bR^k)]
\end{equation}
 and $h_A:\bR^k \to \bR^k$ is given by $h_A(x)=Ax$.
To justify the application of this proposition, we
note that $E\|A\|^{\alpha+\gamma}\leq E\|\Psi\|_{\infty}^{\alpha+\gamma}<\infty$ and $E\|A\|^{\alpha+\gamma}>0$, where $\|A\|=\max_{i \leq k}|\Psi(t_i)|$. (If $E\|A\|^{\alpha+\gamma}=0$ then $P(\Psi(t_i)=0)=1$ for all $i\leq k$, which contradicts \eqref{nonzero-psi}.)

We now prove that the measure $\nu^X$ given by \eqref{def-nux-prod} has marginal measures $\nu_{t_1,\ldots,t_k}^X$. For any $B \in \cB(\barR_0^k)$, we have
\begin{eqnarray*}
\nu_{t_1,\ldots,t_k}^{X}(B)&=&E[\nu_{t_1, \ldots,t_k}\{x \in \bR^k;Ax \in B \cap \bR^k\}] \\
&=& E[\nu_{t_1,\ldots,t_k}\{x\in \bR^k; (\Psi(t_1)x_1,\ldots,\Psi(t_k)x_k) \in B \cap \bR^k\}] \\
&=& E[\nu\{x \in \bD; (\Psi(t_1)x(t_1),\ldots,\Psi(t_k)x(t_k)) \in B \cap \bR^k\}]\\
&=& E[\nu\{x \in \bD; h_{\Psi}(x) \in \{y \in \bD; (y(t_1),\ldots,y(t_k)) \in B \cap \bR^k \}\}] \\
&=& \nu^{X}\{y \in \bD; (y(t_1),\ldots,y(t_k)) \in B \cap \bR^k\}  \\
&=& \nu^{X}(\pi_{t_1,\ldots,t_k}^{-1}(B \cap \bR^k))
\end{eqnarray*}
using \eqref{fin-dim-nu} for the third equality and \eqref{def-nux-prod} for the second last equality.

Next we show that $X$ satisfies condition (b). We only prove (C1). Conditions (C2) and (C3) can be proved similarly. Let $\varepsilon>0$ be arbitrary. If $w''(X,\delta)> a_n \varepsilon$ then there exist some $t_1\leq t \leq t_2$ with $t_2-t_1 \leq \delta$ such that
$|X(t)-X(t_1)|>a_n\varepsilon$ and $|X(t_2)-X(t)|>a_n \varepsilon$. Since
\begin{eqnarray*}
|X(t)-X(t_1)| & \leq & |Z(t)||\Psi(t)-\Psi(t_1)|+|\Psi(t_1)||Z(t)-Z(t_1)| \\
& \leq & \|Z\|_{\infty} |\Psi(t)-\Psi(t_1)|+\|\Psi\|_{\infty}|Z(t)-Z(t_1)|,
\end{eqnarray*}
it follows that
$\|Z\|_{\infty} |\Psi(t)-\Psi(t_1)|>a_n \varepsilon/2$ or $\|\Psi\|_{\infty}|Z(t)-Z(t_1)|> a_n \varepsilon/2$. Similarly,
$\|Z\|_{\infty} |\Psi(t_2)-\Psi(t)|>a_n \varepsilon/2$ or $\|\Psi\|_{\infty}|Z(t_2)-Z(t)|> a_n \varepsilon/2$. Hence
\begin{eqnarray*}
n P(w''(X,\delta)>a_n \varepsilon) & \leq & n P(\|Z\|_{\infty} w''(\Psi,\delta)>a_n \varepsilon/2)+n P(\|\Psi\|_{\infty} w''(Z,\delta)>a_n \varepsilon/2)  \\
& & +2nP(\|Z\|_{\infty}w(\Psi,\delta)>a_n \varepsilon/2)\\
&=:& P_{n,1}(\delta)+P_{n,2}(\delta)+P_{n,3}(\delta).
\end{eqnarray*}
We treat separately the three terms.

For the first term, we note that for any $\theta>0$,
$$P_{n,1}(\delta)
 \leq
n P(\|Z\|_{\infty} \theta>a_n \varepsilon/2)+ n P(\|Z\|_{\infty} w''(\Psi,\delta)1_{\{ w''(\Psi,\delta) >\theta\}}>a_n \varepsilon/2).$$
We take the limit as $n\to \infty$. Using \eqref{tail-conv} for the first term and Lemma \ref{breiman-lemma-R} for the second term, we obtain that, for any $\theta>0$,
$$\limsup_{n \to \infty}P_{n,1}(\delta) \leq c (\varepsilon/2)^{-\alpha}\theta^{\alpha}+c(\varepsilon/2)^{-\alpha}
E[w''(\Psi,\delta)^{\alpha}1_{\{w''(\Psi,\delta)>\theta\}}].$$
Taking $\theta \to 0$, we obtain that $\limsup_{n \to \infty}P_{n,1}(\delta) \leq c(\varepsilon/2)^{-\alpha}
E[w''(\Psi,\delta)^{\alpha}]$.
Take the limit as $\delta \to 0$. By the dominated convergence theorem, $\lim_{\delta \to 0}E[w''(\Psi,\delta)^{\alpha}]=0$, since $\lim_{\delta \to 0}w''(\Psi,\delta)=0$ and $w''(\Psi,\delta) \leq 2 \|\Psi\|_{\infty}$. Hence
$$\lim_{\delta \to 0} \limsup_{n \to \infty}P_{n,1}(\delta)=0.$$

For the second term, we denote by $P_{\Psi}$ the law of $\Psi$ on $\bD$. Since $Z$ and $\Psi$ are independent, we have:
$$P_{n,2}(\delta)=\int_{\bD} n P(\|\psi\|_{\infty} w''(Z,\delta)>a_n \varepsilon/2)  P_{\Psi}(d\psi).$$
Using Lemma \ref{pratt-lemma} (Appendix A), we infer that:
\begin{equation}
\label{ineq-prod}
\limsup_{n \to \infty} P_{n,2}(\delta) \leq \int_{\bD} \limsup_{n \to \infty}
n P(\|\psi\|_{\infty} w''(Z,\delta)>a_n \varepsilon/2)  P_{\Psi}(d\psi).
\end{equation}
To justify the application of this lemma, we note that $$f_n(\psi):=n P(\|\psi\|_{\infty} w''(Z,\delta)>a_n \varepsilon/2)  \leq g_n(\psi):=nP(\|\psi\|_{\infty}\|Z\|_{\infty}>a_n \varepsilon/4),$$
$g_n(\psi) \to g(\psi):=c (\varepsilon/4)^{-\alpha}\|\psi\|_{\infty}^{\alpha}$ (due to \eqref{tail-conv}) and
$$\int_{\bD}g_n(\psi)P_{\Psi}(d\psi)=nP(\|\Psi\|_{\infty} \|Z\|_{\infty}>a_n \varepsilon/4) \to c(\varepsilon/4)^{-\alpha}E\|\Psi\|_{\infty}^{\alpha}=\int_{\bD}g(\psi)P_{\Psi}(d\psi),$$
by Lemma \ref{breiman-lemma-R}. (Note that $\|Z\|_{\infty}$ is regularly varying, and \eqref{nonzero-psi} implies that $P(\|\Psi\|_{\infty}>0)>0$, which forces $E\|\Psi\|_{\infty}^{\alpha+\gamma}>0$.)

We now take the limit as $\delta \to 0$ in \eqref{ineq-prod}.  We apply again Lemma \ref{pratt-lemma} to interchange the limit with the integral. (Both terms are increasing functions of $\delta$, so the limit as $\delta \to 0$ exists.) 
Since $Z$ is regularly varying, $\lim_{\delta \to 0} \limsup_{n}n P(\|\psi\|_{\infty} w''(Z,\delta)>a_n \varepsilon/2)=0$ for any $\psi \in \bD$, and hence
$$\lim_{\delta \to 0} \limsup_{n \to \infty}P_{n,2}(\delta)=0.$$
This second application of Lemma \ref{pratt-lemma} is justified by the fact that $F_{\delta}(\psi):=\limsup_n n P(\|\psi\|_{\infty} w''(Z,\delta)>a_n \varepsilon/2) \leq \limsup_n n P(\|\psi\|_{\infty} \|Z\|_{\infty}>a_n \varepsilon/4):=G(\psi)$ which does not depend on $\delta$.

It remains to treat the third term. Since $\Psi$ has continuous sample paths,
$\lim_{\delta \to 0} w(\Psi,\delta)=0$, and this term is treated exactly as the first term.
$\Box$

\vspace{3mm}

We now consider a finite sum of product terms as in the previous lemma.

\begin{lemma}
\label{sum-lemma}
Let $Z_j=\{Z_j(t)\}_{t \in [0,1]},j=1,\ldots,m$ be i.i.d. processes with sample paths in $\bD$ such that $Z_1 \in {\rm RV}(\{a_n\},\nu,\barD_0)$ and $\alpha>0$ be the index of $Z_1$. Let $\Psi_j=\{\Psi_j(t)\}_{t \in [0,1]},j=1, \ldots,m$ be some processes with continuous sample paths such that $\Psi_1$ is independent of $Z_1$ and
$(\Psi_1, \ldots,\Psi_{j},Z_1, \ldots,Z_{j-1})$ is independent of $Z_j$  for any $j=2,\ldots,m$.
Suppose that there exists a set $T_1 \subset [0,1]$ containing $0$ and $1$ with $[0,1]\verb2\2 T_1$ countable, such that \eqref{nonzero-psi-sum} holds,
and there exists $\gamma >0$ such that $E\|\Psi_j\|_{\infty}^{\alpha+\gamma}<\infty$ for all $j=1, \ldots,m$. Then $X=\sum_{j=1}^{m}\Psi_j Z_j \in {\rm RV}(\{a_n\},\nu^X,\barD_0)$ where
\begin{equation}
\label{def-nux-sum}
\nu^{X}=E[\sum_{j=1}^{m}\nu \circ h_{\Psi_j}^{-1}(\cdot \cap \bD)].
\end{equation}
\end{lemma}

\noindent {\bf Proof:}
Let $T \subset [0,1]$ and $\{\nu_{t_1, \ldots,t_k};t_1, \ldots,t_k \in T,k \geq 1\}$ be the set and the marginal measures given by Theorem \ref{HL-theorem}.(ii) for $Z_1$. We show that $X$ satisfies the conditions of Theorem \ref{HL-theorem}.(ii) with $T_X=T \cap T_1$ instead of $T$ and the measures $\nu_{t_1, \ldots,t_k}^X$ (defined by \eqref{marginals-nux-sum} below) instead of $\nu_{t_1,\ldots,t_k}$.

First we show that $X$ satisfies condition (a). For this, let $t_1, \ldots,t_k \in T_X$ be arbitrary. Note that $(X(t_1),\ldots,X(t_k))^T=\sum_{j=1}^{m}A_j Y_j$ where $A_j$ is the diagonal matrix with entries $\Psi_j(t_1), \ldots, \Psi_j(t_k)$ and $Y_j=(Z_j(t_1), \ldots,Z_j(t_k))^T$. By Lemma 4.3 of \cite{hult-samorod08}, $(X(t_1), \ldots,X(t_k)) \in {\rm RV}(\{a_n\},\nu_{t_1,\ldots,t_k}^X,\barR_0^k)$ where
\begin{equation}
\label{marginals-nux-sum}
\nu_{t_1, \ldots,t_k}^{X}=E[\sum_{j=1}^{m}\nu \circ h_{A_j}^{-1}(\cdot \cap \bR^k)].
\end{equation}
To justify the application of this lemma, we note that $E\|A_j\|^{\alpha+\gamma} \leq E\|\Psi_j\|_{\infty}^{\alpha+\gamma}\linebreak <\infty$ for any $j$, and $P(\cup_{j=1}^{m}\{\|A_j\|>0\})>0$, where $\|A_j\|=\max_{i \leq k}|\Psi_j(t_j)|$. (If $P(\cap_{j=1}^{m}\{\|A_j\|=0\})=1$ then $P(\cap_{j=1}^{m}\{\Psi_j(t_i)=0\})=1$ for any $i=1,\ldots,k$, which contradicts \eqref{nonzero-psi-sum}.) The fact that the measure $\nu^X$ given by \eqref{def-nux-sum} has the marginal measures $\nu_{t_1,\ldots,t_k}^X$
follows as in the case $m=1$ (see the proof of Lemma \ref{product-lemma}). We omit the details.

Next we show that $X$ satisfies condition (b). We only prove (C1). Conditions (C2) and (C3) can be proved by similar methods.

To simplify the notation, we assume that $m=2$. The general result can be proved similarly. Let $\varepsilon>0$ and $\theta>0$ be arbitrary. As in the proof of Lemma 5.1 of \cite{davis-mikosch08}, we use the decomposition:
\begin{eqnarray*}
n P(w''(X,\delta)>a_n \varepsilon) & = & n P(w''(X,\delta)>a_n \varepsilon,\|\Psi_1 Z_1\|_{\infty}>a_n \theta,\|\Psi_2 Z_2 \|_{\infty}>a_n \theta)\\
& + & n P(w''(X,\delta)>a_n \varepsilon,\|\Psi_1 Z_1\|_{\infty}>a_n \theta,\|\Psi_2 Z_2 \|_{\infty} \leq a_n \theta) \\
&+&n P(w''(X,\delta)>a_n \varepsilon,\|\Psi_1 Z_1\|_{\infty} \leq a_n \theta,\|\Psi_2 Z_2 \|_{\infty} > a_n \theta) \\
&+& n P(w''(X,\delta)>a_n \varepsilon,\|\Psi_1 Z_1\|_{\infty}\leq a_n \theta,\|\Psi_2 Z_2 \|_{\infty} \leq a_n \theta)\\
&:=& P_{n,1}(\delta)+ P_{n,2}(\delta)+ P_{n,3}(\delta)+ P_{n,4}(\delta).
\end{eqnarray*}

We treat separately the four terms. For the first term we use the fact that $\|xy\|_{\infty} \leq \|x\|_{\infty}\|y\|_{\infty}$ for all $x,y \in \bD$. We denote by $P_{\Psi_1,\Psi_2,Z_1}$ the law of $(\Psi_1,\Psi_2,Z_1)$. Using the independence between $(\Psi_1,\Psi_2,Z_1)$ and $Z_2$, we have:
\begin{eqnarray*}
P_{n,1}(\delta) &\leq & n P(\|\Psi_1\|_{\infty}\|Z_1\|_{\infty}>a_n\theta, \|\Psi_2\|_{\infty}\|Z_2\|_{\infty}>a_n \theta)\\
&=& \int_{\bD^3} f_n(\psi_1,\psi_2,z_1) dP_{\Psi_1,\Psi_2,Z_1}(\psi_1,\psi_2,z_1).
\end{eqnarray*}
where $f_n(\psi_1,\psi_2,z_1)=n1_{\{\|\psi_1\|_{\infty}\|z_1\|_{\infty}>a_n \theta\}} P(\|\psi_2\|_{\infty}\|Z_2\|_{\infty}>a_n \theta) \to 0$.
By Lemma \ref{pratt-lemma} (Appendix A), it follows that for any $\delta>0$
$$\limsup_{n \to \infty}P_{n,1}(\delta) \leq \int_{\bD^3} \limsup_{n \to \infty}
f_n(\psi_1,\psi_2,z_1) dP_{\Psi_1,\Psi_2,Z_1}(\psi_1,\psi_2,z_1)=0.$$
To justify the application of this lemma, we note that $f_n\leq g_n$ where
$g_n(\psi_1,\psi_2,z_1)=nP(\|\psi_2\|_{\infty}\|Z_2\|_{\infty}>a_n\theta) \to g(\psi_1,\psi_2,z_1)=c\theta^{-\alpha}\|\psi_2\|_{\infty}^{\alpha}$, and
$$\int_{\bD^3} g_{n}dP_{\Psi_1,\Psi_2,Z_1}=nP(\|\Psi_2\|_{\infty}\|Z_2\|_{\infty}>a_n \theta) \to c\theta^{-\alpha}E\|\Psi_2\|^{\alpha}=\int_{\bD^3}g dP_{\Psi_1,\Psi_2,Z_1}.$$
(The last convergence follows by Lemma \ref{breiman-lemma-R} if $P(\|\Psi_2\|_{\infty}>0)>0$, and holds trivially if $\|\Psi_2\|_{\infty}=0$ a.s.)

For the second term, we use the fact that $w''(x+y,\delta) \leq w''(x,\delta)+2\|y\|_{\infty}$ for any $x,y \in \bD$. Hence,
\begin{eqnarray*}
P_{n,2}(\delta) & \leq & n P(w''(\Psi_1 Z_1,\delta)+2\|\Psi_2 Z_2 \|_{\infty}>a_n\varepsilon,\|\Psi_2 Z_2\|_{\infty} \leq a_n \theta) \\
& \leq & n P(w''(\Psi_1 Z_1,\delta)>a_n(\varepsilon-2\theta)).
\end{eqnarray*}
Therefore, if $\theta<\varepsilon/2$, then by Lemma \ref{product-lemma},
$$\lim_{\delta \to 0}\limsup_{n \to \infty}P_{n,2}(\delta) \leq \lim_{\delta \to 0} \limsup_{n \to \infty}nP(w''(\Psi_1 Z_1,\delta)>a_n(\varepsilon-2\theta))=0.$$

The third term is similar to the second term. For the fourth term, we use the fact that $w''(x+y, \delta) \leq 2\|x+y\|_{\infty} \leq 2(\|x\|_{\infty}+\|y\|_{\infty})$ and hence,
$$P_{n,4}(\delta) \leq n P(\|\Psi_1 Z_1\|_{\infty}+\|\Psi_2 Z_2\|_{\infty}>a_n \varepsilon/2,\|\Psi_1 Z_1\|_{\infty} \leq a_n\theta,\|\Psi_2 Z_2\|_{\infty} \leq a_n \theta).$$
The last probability is 0 if $\theta<\varepsilon/4$. The conclusion follows. $\Box$

\pagebreak

\noindent {\bf Proof of Theorem \ref{RV-infinite-sum}:} {\em Step 1.} We show that with probability 1, the series $X(t)$ converges for any $t \in [0,1]$, and the process $X=\{X(t)\}_{t \in [0,1]}$ has sample paths in $\bD$. By applying Theorem 3.1 of \cite{hult-samorod08} to the regularly varying random variables $\{\|Z_j\|_{\infty}\}_{j \geq 1}$, the random coefficients $\{\|\Psi_j\|_{\infty}\}_{j \geq 1}$, and the filtration $\{\cF_j\}_{j \geq 1}$ given by: $\cF_1=\sigma(\Psi_1)$ and $\cF_j=\sigma(\Psi_1,\ldots,\Psi_j,Z_1, \ldots,Z_{j-1})$ for $j \geq 2$, we infer that $\sum_{j \geq 1}\|\Psi_{j}\|_{\infty} \|Z_j\|_{\infty}<\infty$ a.s. Hence, with probability $1$, for any $t \in [0,1]$,
$$|X(t)| \leq \sum_{j \geq 1}|\Psi_j(t)||Z_j(t)| \leq \sum_{j \geq 1}\|\Psi_{j}\|_{\infty} \|Z_j\|_{\infty}<\infty.$$
Let $X^{(m)}=\sum_{j=1}^{m}\Psi_jZ_j$. For any $t \in [0,1]$, $$|(X^{(m)}-X)(t)| \leq \sum_{j \geq m+1}|\Psi_j(t)Z_j(t)| \leq \sum_{j \geq m+1}\|\Psi_j\|_{\infty}\|Z_j\|_{\infty}$$
and hence $\|X^{(m)}-X\|_{\infty} \leq \sum_{j \geq m+1}\|\Psi_j\|_{\infty}\|Z_j\|_{\infty} \to 0$ a.s. Since the uniform limit of a sequence of functions in $\bD$ is in $\bD$, $X \in \bD$ a.s. Since uniform convergence implies $J_1$-convergence, $d_{0}(X^{(m)},X) \to 0$ a.s.

{\em Step 2.} We show that $X \in {\rm RV}(\{a_n\},\nu^X,\barD_0)$, i.e.
$nP(T(a_n^{-1}X) \in \cdot) \stackrel{\hat{w}}{\to}\mu^{X}$ in $\barD_0$, where $\mu^{X}=\nu^{X} \circ T^{-1}$ on $(0,\infty) \times \bS_{\bD}$ and $\mu^X(\{\infty\} \times \bS_{\bD})=0$. By Proposition A2.6.II of \cite{D-VJ}, this is equivalent to showing that:
$$\lim_{n \to \infty}nE[f(T(a_n^{-1}X))] = \int_{\barD_0} f(u) \mu^X(du),$$
for any bounded continuous function $f:\barD_0 \to \bR$ which vanishes outside a bounded set.
Let $f$ be such a function. Suppose that $f$ vanishes outside a set $(r,\infty] \times \bS_{\bD}$ for some $r>0$, and $|f(u)| \leq K$ for all $u \in \barD_0$.

By Lemma \ref{sum-lemma}, we know that for $m$ large enough (for which \eqref{nonzero-psi-sum} holds),
$$\lim_{n \to \infty}nE[f(T(a_n^{-1}X^{(m)}))] =\int_{\barD_0} f(u) \mu^{(m)}(du)$$
where $\mu^{(m)}=\nu^{(m)} \circ T^{-1}$ on $(0,\infty) \times \bS_{\bD}$, with $\nu^{(m)}(\cdot)=E[\sum_{j=1}^{m}\nu \circ h_{\Psi_j}^{-1}(\cdot)]$, and $\mu^{(m)}(\{\infty\} \times \bS_{\bD})=0$.
We claim that
$$\lim_{m \to \infty}\int_{\barD_0} f(u) \mu^{(m)}(du)=\int_{\barD_0} f(u) \mu^X(du).$$
To see this, note that for any bounded measurable function $g:\bD_0 \to \bR$,
$$\int_{\bD_0} g(x)\nu^{X}(dx)=\sum_{j \geq 1}\int_{\Omega}\int_{\bD_0} g(\Psi_j(\omega)y)\nu(dy)P(d\omega).$$
Hence
$$\int_{\barD_0}f(u)\mu^{X}(du)=
\int_{\bD_0}f(T(x))\nu^{X}(dx)=\sum_{j \geq 1}\int_{\Omega}\int_{\bD_0} f(T(\Psi_j(\omega)y))\nu(dy)P(d\omega).$$
Similarly,
$$\int_{\barD_0}f(u)\mu^{(m)}(du)=
\int_{\bD_0}f(T(x))\nu^{(m)}(dx)=\sum_{j=1}^{m}\int_{\Omega}\int_{\bD_0} f(T(\Psi_j(\omega)y))\nu(dy)P(d\omega),$$
and therefore,
\begin{eqnarray*}
& & \left|\int_{\barD_0} f(u) \mu^{(m)}(du)-\int_{\barD_0} f(u) \mu^X(du)\right|\leq  \sum_{j \geq m+1} \int_{\Omega} \int_{\bD_0} |f(T(\Psi_j(\omega) y))|\nu(dy) P(d\omega) \\
& & \leq K \sum_{j \geq m+1} \int_{\Omega}\nu(\{x \in \bD_0; \|\Psi_j(\omega)\|_{\infty}\|x\|_{\infty}>r\}) P(d\omega) \\
& & = K \sum_{j \geq m+1} \int_{\Omega} c\left(\frac{r}{\|\Psi_j(\omega)\|_{\infty}}\right)^{-\alpha}
1_{\{\|\Psi_j(\omega)\|_{\infty}>0\}}P(d\omega)\\
& & = K c r^{-\alpha}\sum_{j \geq m+1}E\|\Psi_j\|_{\infty}^{\alpha} \to 0 \quad \mbox{as} \ m \to \infty,
\end{eqnarray*}
using \eqref{nu-of-V} for the first equality above.
It remains to prove that:
$$\lim_{m \to \infty}\limsup_{n \to \infty} n E|f(T(a_n^{-1}X))-f(T(a_n^{-1}X^{(m)}))|=0.$$
This can be proved similarly to (5.3) of \cite{davis-mikosch08}, using Lemma \ref{sum-lemma} above and Theorem 3.1 of \cite{hult-samorod08}. We omit the details. $\Box$

\appendix

\section{Interchanging limsup with an integral}

Pratt's lemma is a useful tool which allows interchanging a limit with an integral (see \cite{pratt60}). In the present article, we need the following version of Pratt's lemma.

\begin{lemma}
\label{pratt-lemma}
Let $(f_n)_{n \geq 1}$ and $(g_n)_{n \geq 1}$ be some measurable functions defined on a measure space $(E,\cE,\mu)$ such that $0 \leq f_n \leq g_n$ for any $n$, $g_n \to g$ and
\begin{equation}
\label{pratt-cond}
\int_E g_n d\mu \to \int_E g  d\mu<\infty.
\end{equation}
Then
$$\limsup_{n \to \infty} \int_E f_n d\mu \leq \int_E \limsup_{n \to \infty}f_nd\mu.$$
\end{lemma}

\noindent {\bf Proof:}  By Fatou's lemma, $\int \liminf_n(g_n-f_n)d\mu \leq \liminf_n \int (g_n-f_n)d\mu$. By \eqref{pratt-cond}, $g_n$ (and $f_n$) are integrable for $n$ large enough, and $\limsup_n f_n$ is also integrable (being bounded by $g$). Hence, the previous inequality becomes:
$$\int_E gd\mu-\int_E \limsup_n f_n d\mu \leq \liminf_n \int_E g_n d\mu-\limsup_n \int_E f_n d\mu.$$
The conclusion follows by \eqref{pratt-cond}.$\Box$

\end{document}